\documentclass[11pt,letterpaper]{article}

\usepackage{amsmath}
\usepackage{amssymb}
\usepackage{bbm}
\usepackage{hyperref}
\usepackage{mathrsfs}
\usepackage{mathtools}
\usepackage{enumitem}

\usepackage{natbib}
\newcommand{\ignore}[1]{}
\newcommand{\nop}[1]{}
\newcommand*{\eg}{\textit{e.g.}}
\newcommand*{\etal}{\textit{et al.~}}
\newcommand*{\ie}{\textit{i.e.}}

\usepackage{pgfplots}
\usepackage{pgfplotstable}
\pgfplotsset{compat=1.13}  
\usepgfplotslibrary{groupplots}
\usepackage{tikz}

\newtheorem{definition}{Definition}[section]

\newenvironment{blockquote}{%
  \par%
  \medskip
  \leftskip=4em\rightskip=2em%
  \noindent\ignorespaces}{%
  \par\medskip}

\author{
  Justus Isaiah Hibshman\\
  University of Notre Dame\\
  \texttt{jhibshma@nd.edu}
  \and
  Tim Weninger\\
  University of Notre Dame\\
  \texttt{tweninger@nd.edu}
}

\title{Higher Order Imprecise Probabilities and Statistical Testing}
\usetikzlibrary{external}
\begin{document}


\maketitle

\begin{abstract}
We generalize standard credal set models for imprecise probabilities to include higher order credal sets -- confidences about confidences. In doing so, we specify how an agent's higher order confidences (credal sets) update upon observing an event. Our model begins to address standard issues with imprecise probability models, like Dilation and Belief Inertia. We conjecture that when higher order credal sets contain all possible probability functions, then in the limiting case the highest order confidences converge to form a uniform distribution over the first order credal set, where we define uniformity in terms of the statistical distance metric (total variation distance). Finite simulation supports the conjecture. We further suggest that this convergence presents the total-variation-uniform distribution as a natural, privileged prior for statistical hypothesis testing.
\end{abstract}


\let\thefootnote\relax\footnote{\textbf{Acknowledgements:} This research is supported by a grant from the US National Science Foundation (\#1652492)}

\section{Introduction}

Discussing \emph{statistical} evidence means discussing probabilities. When someone uses the term ``probability'' they may mean different things. Some use the idea to refer to a frequency of occurrence. Sometimes this frequency is considered to be inherent to each occurrence of the underlying phenomenon (\eg, some interpretations of the collapse of quantum wave functions). Other times, this frequency is considered to be a more general pattern about trends across repeated occurrences (\eg, repeatedly rolling a die). Others understand probabilities in terms of a measure of a person's confidence that an event would/will/did happen or that something is true/false. For example, the whole field of Bayesian epistemology centers on this latter interpretation. In this paper we focus on this latter use: probabilities interpreted as confidences.

Often, a person wishing to employ mathematical notions of confidence does not have clearly defined confidences (\ie, probabilities) for certain events. Single hypothesis testing is a classic example of this phenomenon. In single hypothesis testing, an experimenter will define a single null hypothesis and obtain an experimental result; then the experimenter will determine whether the experimental result counts as evidence against for or against the null hypothesis. However, to truly determine whether the result is evidence against the null, the experimenter needs to consider the null hypothesis vis-a-vis the alternative hypothesis (\ie the hypothesis that the null hypothesis is false). Richard Royall argues cogently for this requirement~\citeyearpar{royall1997statistical}. Unfortunately, while the null hypothesis (usually) defines a probability function over the experimental result space, the alternative hypothesis is usually too vague to actually define probabilities about experimental results. After all, what is, ``The probability of flipping a coin 20 times and obtaining 7 heads given that the coin flipping is \emph{not} fair?'' Such a question has no answer, but \emph{some kind of answer} is needed in order to meaningfully compare the null hypothesis to the alternative.

One answer is to simply avoid or sidestep situations where probability functions are undefined. In the realm of single hypothesis testing, p-values are such an attempt, but this causes p-values to fail to be a valid measure of evidence. To see excellent exposition on why p-values fail to be a proper formulation of statistical evidence, we refer the reader to work by Richard Royall~\citeyearpar{royall1997statistical} and Wagenmakers~\citeyearpar{wagenmakers2007practical}.

Rather than sidestepping situations where probability functions (\ie, confidences) are undefined, others have suggested using certain agnostic confidences which are intended to express complete ambivalence; these are called uninformative priors. For example, the max entropy prior and the Jeffereys prior are considered to be potentially good options for uninformative priors~\citep{jeffreys1998theory, jaynes1957information, jaynes1982rationale, price2002uninformative}. Unfortunately, the term uninformative prior is a bit of a misnomer, as use of an uninformative prior can seriously bias one's analysis. Additional reading about problems with the use of uninformative priors is provided by Lemoine \citeyearpar{lemoine2019moving}.

While we think the idea of an uninformative prior is suspect, we \emph{do} believe there could be a sort of privileged or natural prior. In this paper, we offer a new model for imprecise confidence updates, and by subsequently analyzing this model, we derive what seems to be a natural prior. Our prior is similar to the Jeffereys prior and the lesser known Hellinger priors~\citep{jeffreys1998theory, shemyakin2014hellinger}. However, it revolves around the use of the total variation statistical metric rather than Fisher information or the Hellinger metric, and provides priors which seem more intuitive when intuition seems relevant. Further, we arrive at our prior by first \emph{ignoring} the idea of using a prior, and then finding that that the prior naturally occurs in our model. Thus, while our discovered prior is not uninformative per-say, it seems like it is natural and therefore to be preferred over other priors.

To develop a model for confidence updates which does not try to sidestep undefined probabilities, we can look to the rich field of Imprecise Probabilities. Models for imprecise probabilities are designed to accommodate the notion that an agent, real \emph{or idealized}, may have vague/imprecise confidences. As Peter Walley eloquently expressed it~\citep{walley1991statistical}:

\begin{blockquote}
When there is little or no information on which to base our conclusions, we cannot expect reasoning (no matter how clever or thorough) to reveal a most probable hypothesis... There are limits to the power of reason.
\end{blockquote}

Imprecise confidences have a rich history and have typically taken two forms. The first form is that of using probability \emph{intervals} rather than point probabilities; the work of Peter Walley is an excellent of example of this kind of model~\citeyearpar{walley1991statistical, walley2000towards}. The second form (which we focus on) is that of using \emph{sets} of probability functions rather than a single probability function. Joyce analogizes such a set to a committee, where different committee members have different confidences~\citeyearpar{joyce2010defense}. Kyburg and Levi produced some of the most famous discussion about sets of probability functions (aka credal sets)~\citep{kyburg1977randomness, kyburg1987bayesian, levi1977direct, levi1975indeterminate}, and others, such as Joyce and Bradley, have subsequently explored and defended the use of credal sets~\citep{joyce2010defense, sep-imprecise-probabilities}.

At present, there is uncertainty concerning how to \emph{update} imprecise confidences after observing an event. When working with probability intervals, the most famous methods are the generalized Bayes' rule~\citep{walley1991statistical}, Dempster's rule~\citep{dempster1967upper, shafer1976mathematical}, and the Geometric rule~\citep{diaconis1978review}. For a recent discussion and comparison of these methods, see Gong and Meng~\citeyearpar{gong2021judicious}. When working with credal sets, the chief method of confidence updates corresponds to what Levi refers to as ``temporal credal conditionalization''~\citeyearpar{levi1978confirmational}. In this paper, we generalize credal set models and then provide a corresponding generalization of temporal credal conditionalization to enable confidence updates with our model.

Current credal set models use a single set of probability functions. Each individual probability function in the credal set can be thought of as a way an agent's confidences \emph{could} be if their confidences were fully precise; the set as a whole represents the agent's imprecise confidences. In the present work, we expand this idea to include credal sets about credal sets (\ie, confidences about confidences). Doing so enables us to address situations where common confidence-update models for credal sets seem to fail. For example, we address the issues of dilation and belief inertia~\citep{seidenfeld1993dilation, sep-imprecise-probabilities}.

Throughout the present work, we seek to provide three levels of depth to accommodate various readers. For those whose primary interest is using the higher order imprecise probability model to perform statistical testing, we provide high-level conceptual descriptions of what our model is doing throughout the paper, as well as an example of how to apply it to a specific case (a binomial test). For those who wish to focus on the formalism and the philosophy, we provide a complete formalism. Finally, for those who wish to understand the formalism but do not have sufficient background, we provide an introduction to background concepts and mathematics in Section~\ref{sec:essential_concepts}.

Our paper has the following major components:
\begin{itemize}[leftmargin=2cm]
    \item[Section~\ref{sec:essential_concepts}] Background
    \item[Section~\ref{sec:existing_issues}] Issues with Existing Credal Set Models
    \item[Section~\ref{sec:our_model}] Higher Order Credal Sets and Confidence Updates
    \item[Section~\ref{sec:conjecture}]  A Formal Conjecture that Almost All the Probability Functions in the Highest Orders of Credal Sets Converge to Form the Same Distribution as Second Order Uniform 
    \item[Section~\ref{sec:hypothesis_testing}] Statistical Testing Using our Model and Conjecture to Provide a Natural Prior Distribution 
    \item[Section~\ref{sec:objections}] Objections and Discussion
\end{itemize}

\section{Background Concepts and Formalisms}\label{sec:essential_concepts}

In this section we briefly cover the conceptual ground needed to understand and apply our mathematical results. In the process, we consider what it means to interpret probabilities as confidences. We introduce existing formalisms for imprecise probabilities, as well as other standard mathematical formalisms we utilize.

\subsection{Probabilities as Confidences}\label{sec:probs_as_confidences}

There are many different ways to interpret the idea of a probability. As we have already discussed, in the present work we consider probabilities as measures of confidence.

When a probability function $P$ from a probability space ($\Omega$, $\mathcal{F}$, $P$) is interpreted as representing confidences, $P$ effectively represents confidences and their relationships \textit{at a single epistemic state}. Confidence about an event $E \in \mathcal{F}$ is expressed as $P(E)$. Relationships between events $A, B \in \mathcal{F}$ are expressed via conditional probabilities, $P(A | B)$, $P(B | A)$.

Once we begin to work with probability functions which are mathematically imprecise in certain ways, we find that two things which are usually treated as conceptually equivalent can sometimes come apart: ``$A$ and $B$ \emph{occur} independently.'' and ``Upon observing $B$, your confidence about $A$ does not change.'' See the work of Couso \etal for examples of the different kinds of independence~\citeyearpar{couso2000survey}. In the present work we do not offer a new definition of conditional probabilities; rather, we address this distinction by considering higher order probabilities. When we consider higher order confidences (confidences about confidences), then the statement ``$A$ and $B$ occur independently.'' can be expressed mathematically in terms of the first-order confidences, and ``Upon observing $B$, your confidence about $A$ does not change.'' can be expressed mathematically in terms of higher order confidences. Without this distinction, the current notion of conditional probabilities is problematic, which we note in Section~\ref{sec:existing_issues}. Our models will address some, but not all of H{\'a}jek's adept argument that the standard ratio definition of conditional probabilities is flawed~\citeyearpar{hajek2003conditional}. We suspect that a definition of conditional probabilities satisfying H{\'a}jek would still work with our model.

\subsection{Intuition for Math of Being (Im)Precise or (Un)Defined}\label{sec:being_undefined}

A confidence that is \emph{imprecise} mathematically parallels the notion of a variable which is \emph{undefined}. 
Many expressions in math are undefined. For instance, the fraction ``$\left(\frac{1}{0}\right)$'' is undefined; the expression ``$\left(\lim_{x \rightarrow 0} \frac{1}{x}\right)$'' is not. There are multiple ways for something to be undefined. In the case of $\frac{1}{0}$, the expression is simply outside the scope of what the division operation is defined upon. Consider the following cases:

\begin{itemize}
    \item If we write, ``Let $x$ be a variable.'' then the value of $x$ is undefined, $x$ could be anything, $2, 5, \{\{1\}, \{2, 3\}, \emptyset \}, \langle 0, 3, 7 \rangle$, purple, etc.
    \item If we write, ``Let $x \in \mathbbm{R}$ be a variable.'' then $x$ is still undefined, but we know a bit more about $x$; we know that $x$ could be $\pi, 1, $ or $-365$ but $x$ cannot be $\{2, 4\}$.
    \item If we write, ``Let $x \in \mathbbm{R}$ where $0 \leq x \leq 100$.'' then once again $x$ is undefined, but we know even more about $x$.
    \item Yet again, if we write, ``Let $x \in \{1, 2, 3, 4, 5\}$.'' then the value of $x$ is undefined, but our possible values for $x$ are very limited.
\end{itemize}

We could provide many, many more examples, but the basic idea should be clear. Interestingly, there are many cases where the value of something can be undefined, yet we still have \emph{some} knowledge about it and can use that knowledge to derive other results. For instance, if we write, ``Let $x \in \{1, 2, 3\}$. Let $y = 4$.'' then we can conclude that $x < y$ even though the value of $x$ is undefined.

In the case of probability formalisms, we can do the same thing. For example, if we are given, ``$(\Omega, \mathcal{F}, P)$ is a probability space where $\Omega = \{0, 1, 2\}$ and $\mathcal{F} = \{\emptyset, \{1\}, \{0, 2\}, \{0, 1, 2\}\}$,'' then we can conclude that $P$ is a function, and while the value of $P(\{1\})$ is undefined, we still know that $0 = P(\emptyset) \leq P(\{1\}) \leq P(\{0, 1, 2\}) = 1$. Similarly, if we are given the additional information that $\cup_{E \in \mathcal{F}} \{P(E)\} = \{0, 0.3, 0.7, 1\}$, then while the value $P(\{1\})$ is still undefined, we know that $P(\{1\}) \in \{0, 0.3, 0.7, 1\}$. To formalize this intuition of being undefined, we can create the notion of a ``values set.''

For any variable or function $x$, we define its values set~$\mathcal{V}(x)$ to be $\mathcal{V}(x) = \{x'\ |\ \left(x' \neq x\right)$ is not defined to be true$\}$. That is, $\mathcal{V}(x)$ excludes all values $x$ is defined not to be and includes all values $x$ could be.

This idea of a values set is mathematically identical to the existing formalisms for imprecise probabilities (introduced in Section~\ref{sec:imprecise_formalism}). In the remainder of the present work, we prefer to speak of confidences being imprecise rather than undefined, both in order to align with the existing literature and because we think that ``imprecise'' corresponds to the notion of a vague mental state better than ``undefined'' does.

\subsection{Probability Spaces, Metrics, Lebesgue Integrals, Etc.}

We now introduce the standard formalisms underlying our results. These formalisms are not necessary to get a high-level overview of how our model works. However, they are needed to understand the full mathematical details. These details are especially relevant to answering the question: ``What is the natural parametrization of a set of probability functions, either explicitly or implicitly?'' 

We begin with the standard definition of a probability space. To do so, we must first define a sigma algebra.

\begin{definition}[Sigma Algebra]
Let $X$ be a set. A \emph{Sigma Algebra} (or $\sigma$-algebra) over $X$ is a set $S \subseteq \mathcal{P}(X)$ meeting the following three conditions:
\begin{enumerate}
    \item $X \in S$.
    \item $S$ is closed under complement: $Y \in S \rightarrow X \setminus Y \in S$.
    \item $S$ is closed under countable union: Let $Y_1, Y_2, Y_3, \ldots \in S$ be a countable sequence. Then $\bigcup_{i} Y_i \in S$.
\end{enumerate}
\end{definition}

\begin{definition}[Probability Space]
Let $X$ be a set and $S$ be a $\sigma$-algebra over $X$. Then the triple ($X$, $S$, $P$) is a probability space if $P$ satisfies the following:
\begin{enumerate}
    \item $P$ is non-negative: $\forall Y \in S.\ P(Y) \geq 0$.
    \item $P(X) = 1$.
    \item $P$ is additive over countable disjoint union: Let $Y_1, Y_2, Y_3, ... \in S$ be any countable sequence where $Y_i \cap Y_j = \emptyset$ for $i \neq j$. Then $P\left(\bigcup_{i} S_i\right) = \sum_i P(Y_i)$.
\end{enumerate}
\end{definition}

When $(X, S, P)$ is a probability space, then we can say that $P$ is a probability measure (or probability function) with respect to $(X, S)$, where $X$ is called the sample space and $S$ is called the event space.

When working with probabilities, it is very useful to be able to sum or integrate over the sample space. This can get complex when the sample space is an infinite set. Integrating over an infinite set requires the following process:

\begin{enumerate}
    \item Begin by defining how far apart elements of $X$ are from each other. Formalize this as a distance metric.
    \item Use the distance metric to generalize the notion of an open interval ($a, b$) to the notion of an open ball and open set.
    \item Define a $\sigma$-algebra $S$ that includes all the open sets. Because $\sigma$-algebras are closed under complement, this means all the closed sets are included as well. When using the $\sigma$-algebra $S$ as the event space, then including all the open sets in $S$ leads to including all the sorts of events that seem intuitively \textit{natural} or related with respect to our distance metric.
    \item Define a notion of the size of elements of $S$. Formalize as a measure over our algebra space.
    \item Use the Lebesgue integral with respect to our measure.
\end{enumerate}

The terms above are defined formally as follows:

\begin{definition}[Distance Metric]
Given a set $X$, a \emph{Distance Metric} on $X$ is a function $d: X \times X \rightarrow \mathbbm{R}^+$ satisfying:
\begin{enumerate}
    \item Identity: $d(a, b) = 0 \leftrightarrow a = b$.
    \item Symmetry: $d(a, b) = d(b, a).$
    \item Triangle inequality: $d(a, c) \leq d(a, b) + d(b, c)$.
\end{enumerate}
\end{definition}

\begin{definition}[Open Ball]
Given a set $X$ with distance metric $d$, an \emph{Open Ball} is a set $Y \subseteq X$ of the form $Y = \{y\ |\ d(x, y) < r\}$ for some $x \in X, r > 0$. Such a ball $Y$ is typically denoted as $B(x, r)$.
\end{definition}

\begin{definition}[Open Set]
Given a set $X$ with distance metric $d$, an open set is a set $Y \subset X$ for which every element of $Y$ is the center of some open ball which itself is a subset of $Y$. Formally, $\forall y \in Y.\ \exists \epsilon > 0.\ B(y, \epsilon) \subset Y.$
\end{definition}

\begin{definition}[Supremum and Infimum]
The infimum (or ``inf'') and supremum (or ``sup'') are effectively min and max functions that can handle some sets which do not have minimums (\eg open sets). In such sets, there might not be a minimum or maximum element, but there might be a \emph{limit} to how large or small the elements get. Given a subset $S$ of a partially ordered set $T$, the \emph{Infimum} of $S$ is the \emph{largest} value $x \in T$ such that $\forall y \in S.\ x \leq y$. Similarly, the \emph{Supremum} of $S$ is the \emph{smallest} value $x \in T$ such that $\forall y \in S.\ y \leq x$. 
\end{definition}

\begin{definition}[Borel Sigma Algebra]
Given a set $X$ with distance metric $d$, the \emph{Borel Sigma Algebra} is the smallest $\sigma$-algebra containing all the open sets in metric space ($X, d$). We indicate the Borel sigma algebra as $\mathscr{B}(X, d)$ or just $\mathscr{B}(X)$ when the metric is obvious.
\end{definition}

\begin{definition}[Measure]
Given a set $X$ with $\sigma$-algebra $S$, a \emph{Measure} for ($X$, $S$) is a function $\mu: S \rightarrow \mathbbm{R}^+ \cup \{\infty\}$ satisfying:
\begin{enumerate}
    \item $\mu(\emptyset) = 0$.
    \item $\mu$ is closed under countable disjoint union: Given a countable sequence $S_1, S_2, S_3, ... \in S$, where $S_j \cap S_k = \emptyset$ for $j \neq k$: $\mu(\bigcup_i S_i) = \sum_i \mu(S_i)$.
\end{enumerate}
\end{definition}

Note that the definition for a measure is very similar to the definition of a probability function. This is because a probability function is simply a measure with 1 as its largest output value.

\begin{definition}[Simple Function]\label{def:simple_function}
Given a set $X$ with $\sigma$-algebra $S$, a function $f: X \rightarrow \mathbbm{R}$ is called a \emph{Simple Function} if $f$ can be expressed as a finite combination of positively-weighted indicator functions for disjoint sets. That is, $f(x) = \sum_{i = 1}^k \alpha_i \mathbbm{1}_{S_i}(x)$ where $S_i \in S$ are disjoint from each other, and where $\mathbbm{1}_{S_i}(x)$ returns 1 if $x \in S_i$ and 0 otherwise.
\end{definition}

\begin{definition}[Measurable Function]
Given spaces ($X_1, S_1$), ($X_2, S_2$) where $S_1$ and $S_2$ are sigma algebras for $X_1$ and $X_2$ respectively, a function $f: X_1 \rightarrow X_2$ is called \emph{Measurable} if: $\forall S \in S_2.\ \{x_1\ |\ f(x_1) \in S\} \in S_1$. 
\end{definition}

Given these definitions, we can now provide a definition of the Lebesgue integral sufficient for present work:

\begin{definition}[Lebesgue Integral]
Given a set $X$ with $\sigma$-algebra $S$, a measure $\mu$ on $S$, and a simple function $f : X \rightarrow \mathbbm{R}^+$ where $f(x) = \sum_{i = 1}^k \alpha_i \mathbbm{1}_{S_i}(x)$ as in Definition~\ref{def:simple_function}, the \emph{Lebesgue Integral} of $f$ on a set $E$ is defined to be:

$$\int_E f\ \text{d}\mu \vcentcolon= \sum_{i = 1}^k \alpha_i\ \mu(S_i \cap E)$$

Similarly, given a measureable function $g : X \rightarrow \mathbbm{R}^+$, where $g$ might not be simple, the \emph{Lebesgue Integral} of $g$ on a set $E \in S$ is defined to be:

$$\int_E g\ \text{d}\mu \vcentcolon= \sup \left\{\int_E h\ \text{d}\mu :\ 0 \leq h \leq g, h \text{ is simple} \right\}$$
\end{definition}

We make use of Lebesgue integration in Section~\ref{sec:our_model} for defining the probability that higher order confidences (\ie, confidences about confidences) implicitly assign to events that the first-order confidences consider.

\subsection{Imprecise Probability Functions - Credal Sets}\label{sec:imprecise_formalism}

Having given an intuition behind the notion of an imprecise probability function,
we now present the existing formalisms for modeling imprecise confidences via credal sets.

\begin{definition}[Credal Set]
Let $\Omega$ be a set and $\mathcal{F}$ be a $\sigma$-algebra over $\Omega$. Then a \emph{Credal Set} is a set of functions $C$ such that $\forall P \in C.\ (\Omega, \mathcal{F}, P)$ is a probability space.
\end{definition}

The above definition admits many interpretations, depending on the context. For example, we can think of $\Omega$ as a set of possible worlds and $\mathcal{F}$ as a collection of events defined over the space of possible worlds. Then a credal set can represent an agent's imprecise confidences about those events. As another example, $\Omega$ can represent all the possible outcomes of a scientific experiment, and $\mathcal{F}$ can represent the events which the experimenter might be able to observe. Then a credal set could represent all the hypotheses the experimenter considers; considering multiple hypotheses can naturally lead to imprecise confidences about what will happen.

We believe the credal set model of imprecise confidences to be extremely general. However, it does not capture \emph{all} the ways one can define a probability function. We explore issues with the model in Section~\ref{sec:existing_issues} and then introduce our expanded model of higher order imprecise probabilities in Section~\ref{sec:our_model}.

A lovely summary of the work on credal sets was penned by Bradley~\citeyearpar{sep-imprecise-probabilities}.

\subsection{Total Variation Distance Metric for Probability Distributions}

Our model requires a metric for distances between probability distributions. We use our metric for two things: Defining Borel $\sigma$-algebras over credal sets, and formulating a notion of a uniform distribution which we use specifically for an important conjecture.

We have chosen to use the total variation distance metric. The total variation metric seems to be the most intuitive possible definition of a distance metric between probability functions, extended to operate on distributions over infinite spaces. Perhaps this is why it is sometimes simply referred to as ``the statistical distance''~\citep{watson2018communication}.

Formally, given probability spaces ($\Omega$, $\mathcal{F}$, $P_A$) and ($\Omega$, $\mathcal{F}$, $P_B$), the total variation distance between $P_A$ and $P_B$ is defined to be:

$$\text{TV}(P_A, P_B) = \sup_{E \in \mathcal{F}} |P_A(E) - P_B(E)|$$

Note that in the above equation, $E$ comes from our event space, which means that it can be an event composed of many smaller events. When the event space $\mathcal{F}$ is finite, then $\text{TV}(P_A, P_B)$ basically corresponds to the overlap of $P_A$ and $P_B$:

\begin{equation*}
\begin{aligned}
\text{When } \mathcal{F} \text{ finite, } \text{TV}(P_A, P_B) &= \text{half the non-overlap between } P_A \text{ and } P_B \\
 &= \frac{1}{2}\left(1 - \sum_{\omega \in \Omega} \min \left\{ P_A(\{\omega\}), P_B(\{\omega\})\right\}\right) \\
 &= \text{half the } L_1 \text{ norm} \\
 &= \frac{1}{2} \sum_{\omega \in \Omega} | P_A(\{\omega\}) - P_B(\{\omega\}) |
\end{aligned}
\end{equation*}

Here we can see that the total variation distance is also a generalization of the $L_1$ distance metric.

\subsection{Summary}

Now that we have formally defined probability functions, credal sets, the statistical distance metric, and integration over infinite sets (\ie, Lebesgue integration), we can use these to formalize how higher order credal sets produce probabilities over lower order spaces. That is to say, given a credal set $\mathcal{C}$ with event space $\mathcal{F}$, a probability function $P$ with $\mathcal{C}$ as its sample space will \emph{implicitly} define a probability function over $\mathcal{F}$ via Lebesgue integration over $\mathcal{C}$. We express the details in Section~\ref{sec:our_model}.

Further, Lebesgue integration and the total variation distance metric will allow us to rigorously define a notion of uniformity over a credal set (see Section~\ref{sec:uniformity}). In turn, this enables us to formulate a fascinating conjecture that in certain conditions the highest order confidences converge to form a single distribution over our original event space (Section~\ref{sec:conjecture}). This convergence gives us an excellent candidate for a natural prior, which we use to formulate statistical tests (Section~\ref{sec:hypothesis_testing}).

Before diving into our model however, we formally describe several of the well known problems with the standard credal set model.

\section{Issues with Existing Credal Set Models}\label{sec:existing_issues}

In this section we discuss two famous issues with how standard credal set models update confidences upon observing events: Dilation and Belief Inertia. We then proceed to our model as an alternative in Section~\ref{sec:our_model}.

\subsection{Dilation}\label{sec:dilation}

Dilation is a process by which standard credal set models unintuitively become \emph{less precise} rather than more precise upon receiving new information \citep{seidenfeld1993dilation, white2010evidential, walley1991statistical, sep-imprecise-probabilities}. Given a credal set $\mathcal{C}$ and an event $E$ we can consider the set of probabilities given to $E$ by $\mathcal{C}$: $\{P(E)\ |\ P \in \mathcal{C}\}$. Now consider the following example:

You have two coins. We denote flipping heads/tails on the first coin with events $H_1$/$T_1$ and heads/tails on the second coin with events $H_2/T_2$. You will toss the first coin such that it has probability $\frac{1}{2}$ of landing heads. The second coin produces heads with some unknown proportion $p \in [0, 1]$. Thus your credal set $\mathcal{C}$ will contain a distinct probability function for each possible value of $p$. We will denote each of these probability functions as $P_p$ for different values of $p$. Expressed formally, we have:

\begin{itemize}
    \item $\forall p \in [0, 1].\ P_p(H_1) = \frac{1}{2}$
    \item $\forall p \in [0, 1].\ P_p(H_2) = p$
\end{itemize}

Now, consider the event $M$ that the second coin toss matches the first toss. Initially, \emph{every} probability function $P_p$ in your credal set gives a probability of $\frac{1}{2}$ to $M$:

$$
\begin{aligned}
P_p(M) &= P_p\left((H_1 \cap H_2) \cup (T_1 \cap T_2)\right) \\
 &= P_p(H_1)P(H_2) + P_p(T_1)P_p(T_2) \\
 &= \frac{1}{2} \cdot p + \frac{1}{2} \cdot (1 - p)\\
 &= \frac{1}{2}
\end{aligned}
$$

Therefore $\{P_p(M)\ |\ P_p \in \mathcal{C}\} = \left\{\frac{1}{2}\right\}$. Now, imagine you observe heads on the first toss (event $H_1$). In the standard credal set formalism, you now obtain a new credal set $\mathcal{C}'$ which contains each probability function from $\mathcal{C}$ conditioned on $H_1$. Consequently, you obtain that $\{{P_p}'(M)\ |\ {P_p}' \in \mathcal{C}'\}$ $= \{P_p(M\ |\ H_1)\ |\ P_p \in \mathcal{C}\} = [0, 1]$. This kind of shift, moving from all probability functions in your credal set giving a single probability such as $\frac{1}{2}$ to a larger range such as $[0, 1]$, is called Dilation.

In this example, it seems like observing a heads on the first toss did not actually change your knowledge about whether the second coin would match the first. However, the space of probabilities you assign to $M$ changed dramatically. Furthermore, if you had observed \emph{tails} rather than heads, you would have gotten the same dilation from $\left\{\frac{1}{2}\right\}$ to $[0, 1]$.

This example makes the original confidences seem too precise or the latter confidences seem too imprecise -- or both.

Our model -- plus a conjecture based on our model -- deals with the issue of dilation in the sense that while dilation can still occur at lower order confidences, the highest order confidences entail first order confidences which effectively do not undergo dilation.

\subsection{Belief Inertia}\label{sec:belief_inertia}

Sometimes, credal set models which are sufficiently imprecise to begin with seem incapable of updating confidences upon receiving evidence. This phenomenon has been termed \textit{Belief Inertia}~\citep{sep-imprecise-probabilities}.

Imagine that you have before you an urn which contains 100 balls painted red, blue, and yellow, but you do not know the proportions. You will repeatedly sample from this urn with replacement. Your initial confidences about which color ball you will sample can be modeled as a credal set containing all valid probability distributions over a mix of 100 red/blue/yellow balls.

Imagine you sample a red ball. Intuitively, this should \emph{in some way} increase your confidence that if you sample again with replacement you will get a red ball. However, under the standard credal set model each \emph{individual} probability function in your credal set models the sampling of balls from the urn as occurring independently. Thus each \emph{individual} probability function does not consider sampling a red ball the second time to be any more likely than it did the first time. Consequently, your credal set for the second sample is \emph{exactly the same} as it was concerning the first sample \emph{and} your credal set for the second sample is \emph{exactly the same} as it would have been if you had sampled a yellow ball rather than a red one. The one caveat is that you might discard the probability function assigning a probability of 0 to red; however this would still be the only belief update, no matter how many times red is then re-sampled.

Here we see a simple and straightforward case where credal set models effectively fail to update at all. Some (\eg, Walley~\citeyearpar{walley1991statistical}) argue that belief inertia is simply a consequence of beginning with totally imprecise confidences and not really a \emph{problem}. Others, such as Rinard~\citeyearpar{rinard2013against}, argue that a completely imprecise credal set does not model a valid mental state in the first place, and thus the credal sets which cause belief inertia should not occur. Still others, such as Vallinder~\citeyearpar{vallinder2018imprecise}, believe that belief inertia is a widespread problem for certain interpretations of imprecise credence models. We view belief inertia as a defect of the standard credal set model, and we offer what we believe to be a better model of imprecise confidences.

Our model addresses the belief inertia issue by introducing higher order confidences about which first order confidences are better models of the sampling process. Consequently, not all probability functions in the new (first order) credal set are treated the same after observing the red ball. Our model - plus a conjecture based on our model - effectively eliminates belief inertia at the highest orders of confidence.

\section{Higher Order Imprecise Probabilities}\label{sec:our_model}

\subsection{The Basics -- Higher Order Credal Sets}\label{sec:our_model_basics}

We offer a general model for imprecise confidences. Let $\Omega$ be a set with $\sigma$-algebra $\mathcal{F}$ as our event space. We define an agent's \textit{first order confidences} to be a standard credal set $\mathcal{C}^1$ where $\forall P^1 \in \mathcal{C}^1.\ (\Omega, \mathcal{F}, P^1)$ is a probability space.

Then, an agent's \emph{Second} Order Confidences are also represented as a credal set $\mathcal{C}^2$, but this time the relevant sample space is $\mathcal{C}^1$, not $\Omega$, and the event space is the Borel $\sigma$-algebra over $\mathcal{C}^1$ with respect to the total variation metric. That is, $\forall P^2 \in \mathcal{C}^2.\ (\mathcal{C}^1, \mathscr{B}(\mathcal{C}^1, \text{TV}), P^2)$ is a probability space.

We can reiterate this definition for Third Order Confidences and so on. In general:

$$\forall i \in \{2, 3, 4, ... \}.\ \forall P^i \in \mathcal{C}^i.\ (\mathcal{C}^{i - 1}, \mathscr{B}(\mathcal{C}^{i - 1}, \text{TV}), P^i) \text{ is a probability space}.$$

We will use $\Omega^i$ to denote $\Omega$ when $i = 1$ and $\mathcal{C}^{i - 1}$ when $i > 1$. Similarly, we use $\mathcal{F}^i$ to denote $\mathcal{F}$ when $i = 1$ and $\mathscr{B}(\mathcal{C}^{i - 1}, \text{TV})$ when $i > 1$.

Any higher order probability function $P^i \in \mathcal{C}^i$ implicitly defines probabilities over the lower order event space(s) as follows. First, for any $j$-th order event $E \in \mathcal{F}^j$, we define the family of measureable functions $p_E^i : \mathcal{C}^i \rightarrow [0, 1]$ for $i \geq j$ where $p_E^i(P^i) = P^i(E)$.

Then we can remember that any $i^\textrm{th}$ order probability function $P^i$ is a measure over our space of $i-1^\textrm{th}$ order confidences. Thus we can use Lebesgue integration to define our higher order confidences as:

$$\forall E \in \mathcal{F}^j.\ \forall P^i \in \mathcal{C}^i \text{ where } i > j.\ P^i(E) \vcentcolon= \int_{\mathcal{C}^{i - 1}} p_E^{i - 1} \ \text{d}P^i$$

Sometimes, we want to express confidences about the combination of events of different orders (their ``intersection''). We can define this in two steps. Step one is almost exactly equivalent to our definition above, just with $\Omega^i = \mathcal{C}^{i - 1}$ replaced with an $i^\textrm{th}$ order event $B \in \mathcal{F}^i$, (which means $B \subseteq \mathcal{C}^{i - 1}$):

$$\forall A \in \mathcal{F}^j.\ \forall B \in \mathcal{F}^i \text{ where } i > j.\ \forall P^i \in \mathcal{C}^i.\ P^i(A, B) \vcentcolon= \int_B p_A^{i - 1} \ \text{d}P^i$$

Then, for the second step, we once again consider a family of functions $q_{A,B}^i : \mathcal{C}^i \rightarrow [0, 1]$ where $q_{A,B}^i(P^i) = P^i(A,B)$. This allows us to define the following:

$$\forall A \in \mathcal{F}^k.\ \forall B \in \mathcal{F}^j.\ \forall P^i \in \mathcal{C}^i \text{ where } i > j \land i > k.\ P^i(A, B) \vcentcolon= \int_{\mathcal{C}^{i - 1}} q_{A,B}^{i - 1} \ \text{d}P^i$$

The above definitions allows us to express conditional confidences in a typical fashion:

$$\forall A \in \mathcal{F}^k.\ \forall B \in \mathcal{F}^j.\ \forall P^i \in \mathcal{C}^i \text{ where } i > j \land i > k.\ P^i(A | B) \vcentcolon= \frac{P^i(A, B)}{P^i(B)}$$

For completeness of notation, we specify that $i^\textrm{th}$ order events do not condition $j^\textrm{th}$ order confidences when $i > j$:

$$\forall A \in \mathcal{F}^k.\ \forall B \in \mathcal{F}^i.\ \forall P^j \in \mathcal{C}^j \text{ where } i > j \land j \geq k.\ P^j(A | B) \vcentcolon= P^j(A)$$

\subsection{Complete Agnosticism}

We say that an agent is completely agnostic about an $i^\textrm{th}$ order credal set if their $(i + 1)^\textrm{th}$ order credal set contains \emph{all} probability functions over their $i^\textrm{th}$ order credal set, their $((i + 1) + 1)^\textrm{th}$ order credal set contains all probability functions over their $(i + 1)^\textrm{th}$ order credal set, and so on.

This notion of complete agnosticism becomes useful for considering what evidence an idealized agent would find in an experimental result, and we return to the idea mathematically in Section~\ref{sec:conjecture}, where we express a conjecture that becomes useful in practice (Section~\ref{sec:hypothesis_testing}). In short, we can model hypothesis testing by letting the first order credal set contain all the possible relevant hypotheses, and letting the higher order confidences be completely agnostic about the hypotheses. Our conjecture can then be interpreted as providing ``whatever statistical evidence the experimental result provides to an agent that is completely agnostic about the hypotheses.''

\subsection{Confidence Updates Upon Observing and Event}\label{sec:confidence_updates}

We now turn to the heart of the higher order probability model, which answers the question: ``Given an idealized agent with credal sets $\mathcal{C}^1, \mathcal{C}^2, \ldots$ and an observed event $E$, what new credal sets ${\mathcal{C}^1}', {\mathcal{C}^2}', \ldots$ should the idealized agent have?''

\subsubsection{The Basic Idea}

As we discussed earlier, the most common model for how a credal set updates upon observing an event $E$ is to take all the original confidence functions and replace them with those same functions conditioned on $E$. As we noted in Section~\ref{sec:existing_issues}, simply performing this kind of update with first-order credences is incomplete. Our basic insight is that, upon observing an event $E$, an agent can perform the standard first order conditioning \emph{and} the agent can perform conditioning with their higher order credal sets. This higher order conditioning effectively amounts to noting which probability functions gave higher confidence to the observed event and prioritizes those probability functions in the future.

A formalism for the higher order conditioning is slightly more complex than standard first order conditioning because the domain of the probability functions changes. Intuitively, for every pair of events $E, F$ from our first order (\ie, standard) event space, when we observe $E$ then the event $F$ maps to a new event $F' = F \cap E$. This kind of mapping is surjective but not injective. That is, multiple distinct events $F$ and $G$ can collapse to a single new event after we observe an event $E$: $F' = F \cap E = G' = G \cap E$.

This collapsing process at the event level means that sometimes multiple distinct probability functions will collapse to a single new probability function once an event has been observed. Because the space of $i\textrm{th}$ order probability functions is itself a sample space for the $(i+1)\textrm{th}$ order confidences, this means that higher order probability functions get a new domain. It is this mapping process that makes the higher order conditioning formalism a bit messier than the standard first order confidence conditioning.

Lastly, we note that while we can \emph{observe} events which our first order confidences consider, we do not normally observe a higher order confidence event. However, we suspect higher order confidence events, such as, ``I just remembered something and now my confidences change.'' might be useful in some contexts. As far as we are aware, our model's formalisms can naturally handle these kinds of observations.

\subsubsection{The Formalism}

Given our original sample spaces $\Omega^i$, $\mathcal{F}^i$, after we observe an event $E \in \mathcal{F}^j$, all event spaces of orders $k < j$ remain the same.

At the order $j$ of the observed event $E$ we update to ${\Omega^j}'$ and ${\mathcal{F}^j}'$ where ${\Omega^j}'$ equals the observed event $E$ and ${\mathcal{F}^j}'$ equals $\{F \cap E\ |\ F \in \mathcal{F}^j\}$. Then:

$${\mathcal{C}^j}' = \{P_E^j\ |\ \exists P^j \in \mathcal{C}^j.\ \forall F' \in {\mathcal{F}^j}'.\ P_E^j(F') = P^j(F'\ |\ E)\}$$

So far, the higher order probability model is more or less what would be expected. When the observed event is a first order event, the first order credal set updates according to standard credal conditioning. Now, for higher order confidences, things get a \emph{bit} more complex, but more or less follow the same pattern. The new sample space for the $i\textrm{th}$ order confidences ($i > j$) is ${\Omega^i}' = {\mathcal{C}^{i - 1}}'$ and the new event space is ${\mathcal{F}^i}' = \mathscr{B}({\mathcal{C}^{i - 1}}', \text{TV})$.




We can now define the way that old probability functions from a credal set map onto new elements when an event $E$ is observed. To do this we will use a family of mapping functions $m_E^k : \Omega^{k + 1} \rightarrow {\Omega^{k +1}}'$.

When $k = j$ (the order at which we observed event $E$), we have that $m_E^j(\omega) = \omega$ when $\omega \in E$ and $m_E^j(\omega)$ is undefined when $\omega \notin E$. When $m_E^j(F)$ is used to express the image of a set $F \in \mathcal{F}^j$, we ignore the undefined outputs and define $m_E^j(F) = F \cap E$.

For $k > j$, we have that $m_E^k$ maps from $\mathcal{C}^k$ to ${\mathcal{C}^k}'$ because $\Omega^{k+1} = \mathcal{C}^k$ and ${\Omega^{k+1}}' = {\mathcal{C}^k}'$. We define $m_E^k$ such that:

$$\forall P^k \in \mathcal{C}^k.\ \forall {F^k}' \in {\mathcal{F}^k}'.\ \left(m_E^k\left(P^k\right)\right)\left({F^k}'\right) = P^k\left({m_E^{(k - 1)}}^{-1}\left({F^k}'\right)\ |\ E\right)$$

This means that the new $P^k$ $\left(\text{\ie, }m_E^k(P^k)\right)$ assigns a new event ${F^k}'$ to the probability of the old event $\left(\text{\ie, the preimage }{m_E^{(k - 1)}}^{-1}({F^k}')\right)$ given the observed event $E$. The new credal set is then:

$${\mathcal{C}^i}' = \{m_E^i(P^i)\ |\ P^i \in \mathcal{C}^i\}$$

With this definition of confidence updates, we get two very desirable properties:

\begin{enumerate}
    \item The $m_E$ mapping preserves conditionals over events:
$$\forall A' \in {\mathcal{F}^k}'.\ \forall B \in \mathcal{F}^j.\ \forall P^i \in \mathcal{C}^i \text{ where } i \geq j, k.\ P^i\left(\left(m_B^k\right)^{-1}\left(A'\right) | B\right) = \left(m_B^i\left(P^i\right)\right)\left(A'\right)$$ 
    \item The $m_E$ mapping preserves the algebra, which is necessary for the higher order conditional to be well defined:
 $${F^k}' \in {\mathcal{F}^k}' \rightarrow \left(m_E^k\right)^{-1}\left({F^k}'\right) \in \mathcal{F}^k$$
    $\left(\text{Remember that for }k > 1,\ \mathcal{F}^k = \mathscr{B}(\mathcal{C}^{k - 1}, \text{TV}) \text{ and } {\mathcal{F}^k}' = \mathscr{B}({\mathcal{C}^{k - 1}}', \text{TV})\right)$
\end{enumerate}

\subsection{Examples: Higher Order Probabilities in Action}

In this section we give some examples of how the higher order probability model works, and we show how it begins to address standard issues with imprecise probabilities like dilation and belief inertia. We believe that \emph{if} a conjecture of ours is true (see Section~\ref{sec:conjecture}) then the higher order probability model \emph{fully} addresses dilation and belief inertia, but \emph{regardless of the conjecture}, we think our model is a step in the right direction, as can be seen in this section.

\subsubsection{Dilation -- Coin Matching}

First, we return to the example given in Section~\ref{sec:dilation}. We have two coins and will toss them with probability $\frac{1}{2}$ of heads for the first coin and a completely unknown probability $p$ of heads for the second coin. We consider the event $M$ that the coin tosses match.

Our sample space is $\Omega = \{H_1H_2, H_1T_2, T_1H_2, T_1T_2\}$ and our event space is $\mathcal{F} = 2^\Omega$. The event $M \in \mathcal{F}$ that the coins match is $M = \{H_1H_2, T_1T_2\}$. We will refer to the event $\{H_1H_2, H_1T_2\}$ simply as $H_1$, the event $\{T_1H_2, T_1T_2\}$ simply as $T_1$, and so on.

Each probability function $P$ in the first order credal set $\mathcal{C}^1$ gives probabilities $\left(a_{P}, b_{P}, c_{P}, d_{P}\right)$ to each of the singleton events $\{H_1H_2\}$, $\{H_1T_2\}$, $\{T_1H_2\}$, $\{T_1T_2\}$ respectively, and the rest of the probabilities are entailed by that assignment. We can distinguish each probability function in our first order credal set via a parameter $p \in [0, 1]$ such that $P_p \in \mathcal{C}^1$ is of the form $\left(\frac{1}{2}p, \frac{1}{2}(1-p), \frac{1}{2}p, \frac{1}{2}(1-p)\right)$; $p$ represents the unknown probability of the second coin toss landing heads according to $P_p$.

As discussed in Section~\ref{sec:dilation}, every first order probability function $P_p \in \mathcal{C}^1$ gives a probability of $\frac{1}{2}$ to the event $M$ before observing $H_1$ or $T_1$. However, after observing $H_1$ \emph{or} $T_1$, the confidence in $M$ becomes $p$ or $(1 - p)$ respectively. In either case, the first order probabilities dilate.

Now, what about the second order probability functions? Originally, each second order probability function also gives a chance of $\frac{1}{2}$ to event $M$. Unfortunately, the second order probabilities dilate, but thankfully, they do not do so to the same degree. If we consider the first order credal set, we find that after observing $H_1$ (or $T_1$) $|\{P\ | P\ \in \mathcal{C}^1 \land P(M\ |\ H_1) = q\}| = 1$ for all $q \in [0, 1]$. However, if we consider our \emph{second} order credal set $\mathcal{C}^2$, we find that things begin to change:

$$q \in \{0, 1\} \rightarrow \left|\left\{P^2\ |\ P^2 \in \mathcal{C}^2 \land P^2(M\ |\ H_1) = q \right\}\right| = 1$$

$$q \in (0, 1) \rightarrow \left|\left\{P^2\ |\ P^2 \in \mathcal{C}^2 \land P^2(M\ |\ H_1) = q \right\}\right| = |\mathbbm{R}|$$

Thus we see that higher order confidences do not dilate \textit{as much} as first order confidences. Less confidence is placed at the extremes (0 and 1) and more is placed toward the middle. This is because a higher order confidence function represents a combination of lower order probability functions, so the probabilities tend to \textit{average out} as the order of confidence gets higher. Only one second order distribution can maintain a confidence of $1$ in $M$ after observing $H_1$: the distribution $P^2$ satisfying $P^2(\{P^1 | P^1 \in \mathcal{C}^1 \land P^1(M\ |\ H_1) = 1\} = \{P_1\}) = 1$. We discuss a formal conjecture about this process of \textit{averaging out} in Section~\ref{sec:conjecture}.

Ideally, we would analyze this lessening of dilation in terms of a uniform measure rather than in terms of cardinality. However, we have yet to be able to calculate a well defined uniform measure over a full second (or higher) order credal set when the first order credal set contains infinitely many probability functions; we conjecture it exists though. See Section~\ref{sec:conjecture} for more information.

We can see from this example that the higher order probability model shows promise in dealing with dilation but does not presently remove dilation altogether. However, if our conjecture proves true, then we find almost all of the \textit{highest} order confidence functions undergo no dilation-like behavior (\eg, almost all $P^i(M\ |\ H_1) = \frac{1}{2}$ in this scenario as $i \rightarrow \infty$).

\subsubsection{Belief Inertia -- Urn Sampling}

We return to the example given in Section~\ref{sec:belief_inertia}. We have an urn containing 100 balls colored red, yellow, and/or blue. We do not know the proportions. We will sample with replacement $n$ times.

The sample space is $\Omega = \{\text{Red}, \text{Yellow}, \text{Blue}\}^n$, and the event space is $\mathcal{F} = 2^\Omega$. We can represent a probability function in our first order credal set via just three numbers, $($Probability of red, Prob of yellow, Prob. of blue$)$. The first order credal set will be:

$$\mathcal{C}^1 = \left\{\left(\frac{R}{100}, \frac{Y}{100}, \frac{B}{100}\right)\ |\ R, Y, B \in \mathbbm{N} \cup \{0\} \land R + Y + B = 100 \right\}$$

The higher order credal sets will be completely agnostic over their respective event spaces.

Let the first sample be a red ball. As with standard credal set models, the first order credal sets effectively do not change, because each first order probability functions models the samples as being independent. However, if we consider the space of second order credal sets, we find that almost all of them consider red balls more likely than they did previously after observing the first red ball (where this \textit{almost all} is in terms of the third order uniform distribution over the second order credal set). 
Furthermore, if our conjecture is true (Section~\ref{sec:conjecture}), then as the order increases, almost all of the confidence functions approach a uniform distribution over the first order credal set, which means they move from a confidence of $\frac{1}{3}$ that a red ball will be drawn to a confidence of $\frac{91}{180} = 0.50\bar{5}$ that a red ball will be drawn. Further, if a yellow is observed next, both red and yellow would get a confidence of $\frac{181}{450} = 0.40\bar{2}$. Or, instead, if a red ball is drawn twice (no yellow), then the confidence in red would go to $\frac{24841}{40950} \approx 0.6067$.

Thus, once again, the higher order probability model begins to address the issues with credal sets and, if our conjecture is true, seems to completely address the issue of belief inertia.

\section{Higher Order Agnostic Confidence Convergence Conjecture}\label{sec:conjecture}

Having established our model and given some examples, we now turn to our conjecture about how the higher order probability model behaves in the limiting case as the order of confidences being considered becomes arbitrarily large.

We observe in \emph{finite} simulations that the distributions which higher order probability functions entail over the first order event space tend to converge. In this section, we formalize what this would mean in the infinite case into a formal conjecture. This conjecture provides the basis for our application to statistical testing.

\subsection{Uniformity with Respect to a Distance Metric}\label{sec:uniformity}

Formally defining a uniform distribution over an infinite set is not quite as easy as it may sound at first. We cannot simply say, ``Make each element equally likely.'' For example, consider sampling uniformly from the interval [0, 1]. If we just say that, ``Each element is equally likely.'' then someone could say, ``Since each element is equally likely and a bijection exists between the interval [0, 0.2] and the interval (0.2, 1], then those intervals have the same number of elements and the chance of sampling from [0, 0.2] is equal to the chance of sampling from (0.2, 1].'' This is obviously problematic, as it both defies our understanding of uniformity and leads to contradictions. Thus, more careful formalisms are needed.

The idea of a uniform distribution is typically formalized and understood in terms of the Lebesgue measure $\mathcal{L}$. That is, sampling \textit{uniformly} is typically understood to mean that the chance of sampling from a subset $S$ of a sample space $X$ is just $\frac{\mathcal{L}(S)}{\mathcal{L}(X)}$. In turn, the Lebesgue measure is defined to be the measure which provides a standard notion of volume, and a standard notion of volume is itself defined in terms of a standard notion of distance: Euclidean distance. Returning to above example, the volume of [0, 0.2] is not equal to the volume of (0.2, 1] because the Euclidean distance between 0 and 0.2 is not the same as the distance between 0.2 and 1.

The Lebesgue measure $\lambda$ is \emph{typically} defined for $\mathbbm{R}^n$ in terms of the radii of open balls. For example, it can be thought of as a measure $\lambda$ satisfying the following requirements:

\begin{enumerate}
    \item For any open ball $B(x, r)$, the measure (\ie, size) $\lambda(B(x, r))$ is finite and $> 0$ if and only if $r$ is finite and $> 0$.
    \item All open balls of the same radius have the same measure (\ie, same size).
\end{enumerate}

In more detail, the Lebesgue measure effectively assigns open balls of radius $r$ in $\mathbbm{R}^n$ a size of $r^n$ and assigns any other set the sum of the sizes of the smallest (countable) collection of open balls that will cover the set. However, while this is the more traditional way to understand the Lebesgue measure, there is \emph{another} way which we find more useful for our context.

We can consider an $m$-dimensional parametrization of the points in a set to be measured. When working in $\mathcal{R}^n$, this parametrization is usually just the $m = n$ coordinates of a point. The Lebesgue measure for this $m$-dimensional object can then be thought of as an integral of the $(m-1)$-dimensional Lebesgue measure over the $(m-1)$ dimensional slices of the $m$-dimensional parametrization. However, the slices are not necessarily equally thick with respect to the parameter $x_m$. Rather, a slice's \textit{thickness} is determined by the rate at which the distance metric changes with respect to the parametrization. The Euclidean distance metric has the advantage that the measure for a set obtained in this way does not change with a different parametrization. We will formalize what we mean by all this when we get to the higher order probability model for uniform distributions over credal sets in Section~\ref{sec:defining_uniformity_over_credal_sets}.

\subsubsection{(Im)Possibility of Uniformity}\label{sec:im_possibility_of_uniformity}

Recall that a straightforward way to define a Lebesgue-like measure is to require that all open balls of finite, non-zero radius $r$ have the same finite non-zero size. Sometimes, a metric space does not allow for a measure to be defined in this way. For instance, consider the Euclidean distance metric $L_2$ on the space of $\mathbbm{R}^\infty$. Assume for the sake of contradiction we have a measure $\mu : \mathscr{B}(\mathbbm{R}^\infty, L_2) \rightarrow \mathbbm{R}$ which meets the two criterion described in the section immediately above. Now consider any open ball $B(x, r) \in \mathscr{B}(\mathbbm{R}^\infty, L_2)$ where the radius $r > 0$ is finite. By our assumption, $\mu(B(x, r))$ is also finite and non-zero. However, \emph{because} we are in $\mathbbm{R}^\infty$, we know that we can fit \emph{countably infinitely} many non-overlapping balls of some radius $0 < r' < r$ into $B(x, r)$. Again, given our original assumption, each of these balls of radius $r' > 0$ has the same, non-zero measure. Because $B(x, r)$ contains all these balls and $\mu$ must be additive for countable union, $B(x, r)$ must have infinite measure. Contradiction. Consequently, $\mathbbm{R}^\infty$ with the Euclidean distance metric cannot have a Lebesgue-like measure defined in terms of ball radii. The Lebesgue measure is undefined for $\mathbbm{R}^\infty$.

Consequently, to define a uniform measure over an infinite set we may need something other than Lebesgue.

\subsubsection{Simple Parametrizations}\label{sec:simple_parametrizations}

Before we turn to the defining uniformity over credal sets (Section~\ref{sec:defining_uniformity_over_credal_sets}), we first define a simple parametrization of a credal set.


Let $\mathcal{C}$ be a credal set and $d$ be a distance metric over $\mathcal{C}$. Let $\mathcal{I}$ be a cross-product of $n$ closed intervals: $\mathcal{I} = [a_1, b_1] \times [a_2, b_2] \times ... \times [a_n, b_n]$.

Then we call a function $p : \mathcal{I} \rightarrow \mathcal{C}$ a \textit{simple parametrization} of $\mathcal{C}$ \emph{if} $p$ is a measurable function with respect to $(\mathcal{I}, \mathscr{B}(\mathcal{I}, L_2)) \rightarrow (\mathcal{C}, \mathscr{B}(\mathcal{C}, d))$ (where $L_2$ is the Euclidean distance metric).

\subsubsection{Uniformity Over Credal Sets}\label{sec:defining_uniformity_over_credal_sets}

Can we define a uniform distribution over a space of probability functions -- that is, over a credal set? In particular, can we do so with the respect to a particular metric $d$?

More specifically, can we define a uniform distribution $U: \mathscr{B}(\mathcal{C}^i, d) \rightarrow [0, 1]$ for an $i^\textrm{th}$ order credal set $\mathcal{C}^i$ where the functions $P \in \mathcal{C}^i$ form a probability space $(\Omega^i, \mathcal{F}^i, P)$ with a sample space $\Omega^i$ of cardinality $|\Omega^i| = |\mathbbm{R}|$? At first glance, the answer might seem to be that we cannot. After all, it might seem like a probability function for an infinite event space $\Omega^i$ is basically an infinite-dimensional vector, and as we already noted in Section~\ref{sec:im_possibility_of_uniformity}, problems arise when we try to follow the standard course for defining a uniform distribution over $\mathbbm{R}^\infty$ when our distance metric is Euclidean.

However, the space of probability functions has some important constraints on it that still allow us to define a uniform measure in another fashion. When the size of our sample space $|\Omega^i| \leq |\mathbbm{R}|$, then $|C^i| \leq |\mathbbm{R}|$ even though the space of all \emph{functions} over $|\Omega^i|$ has cardinality $|2^\mathbbm{R}|$. This can be seen by the fact that each probability distribution over a set of cardinality $|\mathbbm{R}|$ is uniquely defined by its cumulative distribution function (cdf), cdf's are right-continuous functions, and the space of all right-continuous functions over $\mathbbm{R}$ has the cardinality $|\mathbbm{R}|$.

This fact gives us two crucial properties: Firstly, if our first order sample space $\Omega^1$ has cardinality $|\Omega^1| \leq |\mathbbm{R}|$, then inductively we get that \emph{all} higher order credal sets have cardinality $|\mathcal{C}^i| \leq |\mathbbm{R}|$. Secondly, this means that we can define a uniform distribution over our first order credal set (and perhaps all higher orders as well). We can formalize uniformity as follows:

Let $\Omega$ be a set with $\sigma$-algebra $\mathcal{F}$ where $|\Omega| \leq |\mathbbm{R}|$ and the event space $\mathcal{F}$ \emph{also} has a cardinality of $|\mathbbm{R}|$ or less. Note that a Borel sigma algebra on a set of size $|\mathbbm{R}|$ \emph{also} has cardinality $|\mathbbm{R}|$~\citep{kannai2019elementary}. $\Omega$ will be our sample space and $\mathcal{F}$ our event space.

Now, consider a credal set $\mathcal{C}$ for $\Omega, \mathcal{F}$. As already discussed, we know that $|\mathcal{C}| \leq |\mathbbm{R}|$. In this paper, we will define uniformity over $\mathcal{C}$ with respect to some distance metric $d$, and will focus on when $d$ is the Total Variation metric.

\textbf{Case 1: $\mathcal{C}$ is finite.}

In this case, simply define our uniform distribution $U_d$ over $\mathcal{C}$ in the standard fashion: Each element of $\mathcal{C}$ has probability $\frac{1}{|\mathcal{C}|}$.

\textbf{Case 2: $|\Omega| \leq |\mathbbm{R}|$ and $|\mathcal{C}| = |\mathbbm{R}|$ with a simple parametrization $p$.}

Because $p$ is a simple parametrization for $\mathcal{C}$, then by definition it has domain and range $p : I_1 \times I_2 \times ... \times I_n \rightarrow \mathcal{C}$ for some $n$. Also, $p$ is continuous with respect to our distance metric $d$. In this case the $I_i$ are closed intervals $[a_i, b_i]$. We represent the lower and upper bounds of interval $I_i$ as $I_i[0]$ and $I_i[1]$ respectively.

We can now define two families of useful functions. First we define our \textit{thickness functions}. Let $t_{d,k} : I_1 \times I_2 \times ... \times I_n \rightarrow \mathbbm{R}$ be defined for $k \in [n]$ such that:

$$t_{p,d,k}(\langle x_1, x_2, \ldots, x_n \rangle) = \lim_{y \rightarrow x_k} \frac{d(p(\langle x_1, \ldots, x_n \rangle), p(\langle x_1, \ldots , x_{k - 1}, y, x_{k + 1}, \ldots , x_n \rangle))}{|x_k - y|}$$

The reader should note that the definition for $t_{p,d,k}$ is essentially the definition of a partial derivative, except that the numerator has a distance metric other than the Euclidean metric. Now we can define our second family of functions, a class of measures, $U_{p,d,k}^*$ for $k$-dimensional sub-intervals of the parameter space:

$$U_{p,d,k}^*: (I_1 \times I_1) \times (I_2 \times I_2) \times ... \times (I_k \times I_k) \times I_{k + 1} \times ... \times I_n \rightarrow \mathbbm{R}$$

Furthermore:

$$U_{p,d,1}^*(\langle J_1, x_2, x_2, ..., x_n \rangle) = \int_{a_1}^{b_1} t_{d,1}(\langle y, x_2, x_3, ..., x_n \rangle)\ \text{d}y$$

\noindent where $J_1 \subseteq I_1$ is an interval $[a_1, b_1]$. This basically means, ``Integrate the rate at which our distance metric changes relative to our parameter over a line in the parameter space.'' This gives us the length of the line in terms of $d$ rather than in terms of our parametrization. We can inductively define this kind of integration for higher orders by specifying for $k > 1$ that:

$$U_{p,d,k}^*(\langle J_1, J_2, \ldots, J_k, x_{k+1}, \ldots, x_n \rangle) = \int_{J_k[0]}^{J_k[1]} \int_{J_1 \times \ldots \times J_{k - 1} ; \langle y, x_{k + 1}, x_{k + 2}, \ldots, x_n \rangle } t_{p,d,k}\ \text{d}U_{p,d,k - 1}^*\ \text{d}y$$

The inner integral basically means, ``For a $(k-1)$-dimensional sub-interval of our parameter space, we can think of this $(k-1)$-dimensional sub-interval as a slice of a $k$-dimensional sub-interval. We could consider just using $U_{p,d,k-1}^*$ to give us a size of the slice, but the slice might not be equally thick in terms of the $k$th dimension. Thus, we will get the size of that slice by integrating the $k$-dimensional thickness function $t_{p,d,k}$ over the slice. The measure $U_{p,d,k-1}^*$ tells us how large sub-parts of the slice are as we integrate.'' Then, the outer integral basically means, ``For a $k$-dimensional sub-interval of our parameter space, we can obtain the size of that sub-interval by integrating the size of all the different $(k-1)$-dimensional slices.''

Given these measures for sub-intervals of the parametrized space, we can now define a measure similar to the Lebesgue measure for any element of our credal set space.

$$U_{p,d}^* : \mathscr{B}(\mathcal{C}, d) \rightarrow \mathbbm{R}$$

\noindent where:


$$
U_{p,d}^*(T) = \inf \left\{\sum_i U_{p,d,n}^*(S_i) : \text{Countable sub-intervals } S_1, S_2, ... \textrm{ of } \mathcal{I} \text{ s.t. } T \subseteq \bigcup_i p(S_i) \right\}
$$

Last but not least we can easily normalize $U_{p,d}^*$ into a probability measure $U_{p,d}$ by specifying that $U_{p,d}(S) = \frac{U_{p,d}^*(S)}{U_{p,d}^*(\mathcal{C})}$. This will be our uniform probability measure.

Note that our \textit{uniform} measure still is denoted in terms of our parametrization $p$; we write $U_{p,d}$ rather than just $U_d$. This is because \emph{for some distance metrics}, $U_{p,d}$ could hypothetically vary depending on $p$. \emph{However}, as far as we are aware, $U_{p,\text{TV}}$ is the same for all simple parametrizations $p$. Thus, when speaking about uniformity with respect to the total variation metric, we will simply write $U_\text{TV}$ with some arbitrary simple parametrization $p$ implied.


\textbf{Other Cases: (\eg, $|\mathcal{C}| = |\mathbbm{N}|$).}

We think other cases are quite interesting, but we do not consider definitions for them in this paper as they are not needed for our main results.

\subsection{Total Variation Uniform is Normal}\label{sec:U_TV_is_normal}

We find it worthwhile to note that in cases where we can have a concrete expectation about what a \textit{uniform} distribution over a first order credal set should look like, $U_\text{TV}$ seems to provide such a distribution. For example, if our first order credal set is for a single coin toss with an unknown proportion $p \in [0, 1]$, $U_\text{TV}$ treats each value for $p$ as equally likely. Similarly, for any distribution over a finite sample space of size $n$, we can think of a first order distribution as a vector of $n$ elements summing to $1$. Our definition of uniformity over this space given by $U_\text{TV}$ produces the same uniformity over this space that the Lebesgue measure produces: our natural instinct for uniformity over a plane/simplex.

\subsection{The Sub-Conjecture}\label{sec:sub_conjecture}

Before proceeding to our main conjecture, we first offer a sub-conjecture. Our sub-conjecture effectively states that a uniform distribution is well defined for higher order credal sets.

Let $\Omega^i$ be an $i$th order sample space with event space $\mathcal{F}^i$ where $(\Omega^i, \mathcal{F}^i)$ is finite or simply parametrizeable. Let our $i^\textrm{th}$ order credal set $\mathcal{C}^i$ contain \emph{all} the possible probability functions for sample/event space $\Omega^i, \mathcal{F}^i$.

\textbf{Minimal Sub-Conjecture: } There exists a simple parametrization for $\mathcal{C}^i$ with respect to the total variation metric.

If this sub-conjecture is true, then assuming we can define a uniform distribution over our first order credal set, we get a uniform distribution for \emph{any} order of credal set.

Our main conjecture assumes this minimal sub-conjecture is true. Before turning to our main conjecture, we offer some small steps toward proving the minimal sub-conjecture via a stronger sub-conjecture.

\textbf{Strong Sub-Conjecture:} There exists a total ordering $\precsim$ on $\mathcal{C}^i$ with min and max elements $P_L, P_U$ such that the following two conditions hold:

\begin{enumerate}
    \item $\forall P_L \precsim P \prec P_U.\ \forall \epsilon > 0.\ \exists P_u \succ P.\ \forall P \precsim Q \precsim P_u.\ \text{TV}(P, Q) < \epsilon.$
    \item $\forall P_L \prec P \precsim P_U.\ \forall \epsilon > 0.\ \exists P_l \prec P.\ \ \forall P_l \precsim Q \precsim P.\ \text{TV}(P, Q) < \epsilon.$
\end{enumerate}
This strong sub-conjecture is basically a continuity requirement for a total ordering. If the strong sub-conjecture holds true, then the minimal sub-conjecture holds true. To see this, consider an ordering satisfying the strong sub-conjecture. This ordering entails the existence of a family of bijections $f: \mathcal{C}^i \rightarrow [0, 1]$ satisfying $P \precsim Q \leftrightarrow f(P) \leq f(Q)$. These bijections are simple parametrizations of $\mathcal{C}^i$.

It is fairly trivial to construct a total ordering on $\mathcal{C}^i$ with max and min elements, but it is much harder to ensure that the ordering is continuous with respect to TV as required by the strong sub-conjecture.

\subsection{The Main Conjecture}\label{sec:conjecture_main_subsection}

Assume the minimal sub-conjecture from Section~\ref{sec:sub_conjecture}.

Let $\mathcal{F}$ be a first order event space. Let an agent's first order credal set $\mathcal{C}^1$ be finite or simply parametrizeable, and let the agent's higher order confidences be completely agnostic about their first order confidences.

By applying the minimal sub-conjecture inductively, we get that \emph{every} credal set is (finite or) simply parametrizeable and thus every credal set has a uniform distribution defined over it. We denote the uniform distribution over the $i^\textrm{th}$ order credal set $\mathcal{C}^i$ with $U_\text{TV}^{i + 1}$.

\textbf{We conjecture that:}

$$\lim_{i \rightarrow \infty} U_\text{TV}^{i+1}\left(\left\{P^i\ |\ P^i \in \mathcal{C}^i \land \exists E \in \mathcal{F}.\ P^i(E) \neq U_\text{TV}^2(E) \right\}\right) = 0$$





\emph{In other words}, when an agent is completely agnostic about their first order confidences, then their confidences converge in the limit of confidence order. As an agent employs higher and higher order confidences, their \textit{highest} order confidences almost all produce the same distribution over the first order event space. Furthermore, they approach a well-defined second order probability function, $U_\text{TV}^2$.

Note that if this conjecture is true we immediately get that it applies to confidence updates as well: When a first order event $E$ is observed, then almost all the highest order confidences in the new credal sets (as measured by $m_E^i(U_\text{TV}^i)$ as $i \rightarrow \infty$) will form the same confidences for any event $F$ as $U_\text{TV}^2(F\ |\ E)$.

\subsection{Finite Simulation Example}\label{sec:convergenc_simulation}

\begin{figure}[t]
\centering
\begin{tabular}{c c}
\input{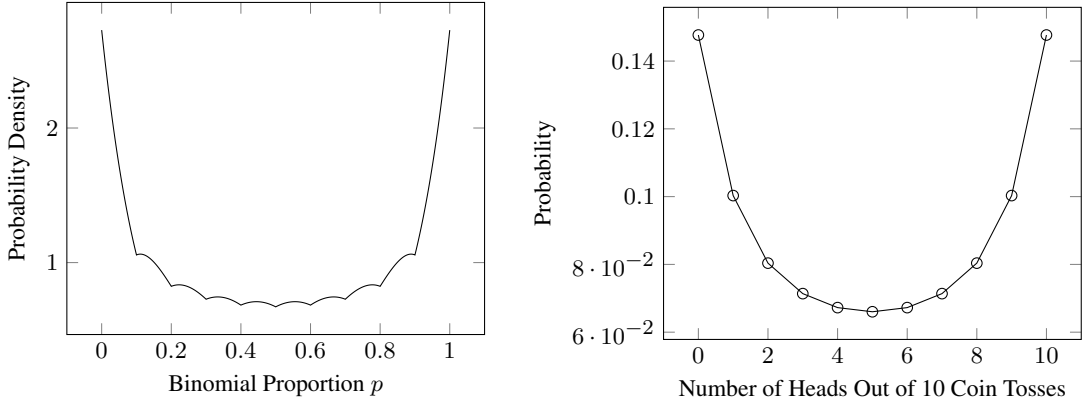} & 
\begin{tikzpicture}
\begin{axis}[
	xlabel=Number of Heads Out of 10 Coin Tosses,
	ylabel=Probability,
	width=0.5\linewidth,height=6cm,
    legend style={at={(0.0,.91)},anchor=west}
    ]

\addplot[color=black, mark=o] coordinates {
	(0, 0.147688334187352)
(1, 0.100306540575736)
(2, 0.0803729874336117)
(3, 0.0713748966559957)
(4, 0.0672460670033027)
(5, 0.0660223482880037)
(6, 0.0672460670033027)
(7, 0.0713748966559957)
(8, 0.0803729874336117)
(9, 0.100306540575736)
(10, 0.147688334187352)

};


\end{axis}
\end{tikzpicture}
\end{tabular}
\caption{The PDF of the distribution that Total-Variation-Uniform ($U_\text{TV}^2$) forms over the binomial distributions for 10 coin tosses (left). The distribution that $U_\text{TV}^2$ implicitly defines over heads from 10 coin tosses (right). Our conjecture in Section~\ref{sec:conjecture_main_subsection} says that higher order probability functions will almost all converge to form the same distribution over heads as this distribution.
}\label{fig:U_TV_2_binomial}
\end{figure}

We can simulate our conjecture in the case that a first order credal set is over a finite sample space.

Recall that in finite cases the total variation metric coincides with the $L_1$ distance metric (TV($P_A, P_B) = \frac{1}{2}L_1(P_A, P_B)$). We use the method of Calafiore \citep{lacko2012conditional, calafiore1998uniform} to randomly generate vectors on an $n$-dimensional $L_1$-simplex in an $L_1$-uniform fashion; this gives a representative sample of the space of probability functions over a sample space of size $n$.

We pick a finite sample space $\Omega$ and a simply parametrizeable first order credal set $\mathcal{C}^1$. We pick an event of interest $E \in \mathcal{F} = 2^\Omega$. Our simulation will run as follows:

\begin{enumerate}
    \item Randomly sample distributions from $\mathcal{C}^1$ via the distribution $U_\text{TV}^2$. This forms our finite, representative approximation of $\mathcal{C}^1$, which we shall refer to as $\bar{\mathcal{C}}^1$.
    \item Randomly sample distributions from the full credal set over $\bar{\mathcal{C}}^i$ in a TV-uniform fashion using Calafiore's method. These samples form $\bar{\mathcal{C}}^{i+1}$.
    \item For $i \in [k]$, plot the space of values $P^i(E)$ for $P^i \in \bar{\mathcal{C}}^i$ as weighted by $\bar{U}_\text{TV}^{i+1}$.
\end{enumerate}

We find that if we make our sample sizes large enough, thereby making the samples more accurate approximations of the infinite case, then at higher orders the probabilities for $E$ tend to converge to $U_\text{TV}^2(E)$, just as our conjecture says they should.

In particular, we consider the set of binomial distributions for 10 coin tosses as our first order credal set. We then look at the probabilities that higher order probability functions give to specific events, and compare those with the probability that $U_\text{TV}^2$ gives to those events. In Figure~\ref{fig:U_TV_2_binomial}, we depict $U_\text{TV}^2$ both over the first order credal set (the binomials) on the left and over the first order event space (the coin tosses) on the right. 

While $U_\text{TV}^2$ bears a resemblance to some depictions of the Jeffereys prior, there are two key differences: The Jeffereys prior's PDF is unbounded as it approaches $p = 0$ and $p = 1$, and the $U_\text{TV}^2$ distribution has sharp points where the proportion $p$ is a fraction $\frac{k}{10}$ for $k \in \{0, 1, ..., 10\}$; these are the proportions that maximize the chance of getting $k$ heads for some $k$.

\begin{figure}[t!]
\centering
    \input{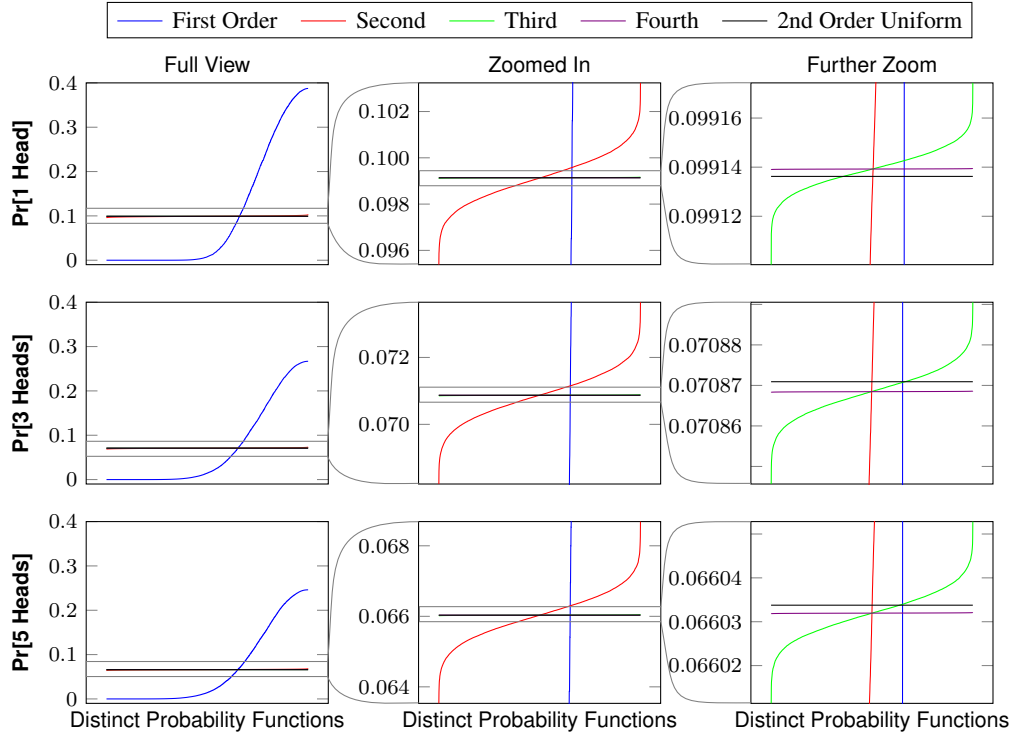}
\caption{Simulating higher order convergence of confidence functions. Here we display the probabilities that different orders of credal sets assign to three distinct possible outcomes from 10 coin tosses: getting 1, 3, and 5 heads, shown in the top, middle, and bottom rows respectively. The x axes contain a representative sample of the $1^\textrm{st}-4^\textrm{th}$ order probability functions, ordered by the probability which they ascribe to the event of that row. Successive columns zoom in on the previous column to show the higher orders more clearly. For each event, as the order of confidence increases, more and more probability functions approach assigning the same probability to the event that second order uniform ($U_\text{TV}^2$) assigns (the black line).}\label{fig:binomial_convergence}
\end{figure}

In Figure~\ref{fig:binomial_convergence} we depict how in finite simulations, higher order credal sets converge on probabilities given to events, such as obtaining 3 heads. These finite simulation results correspond with our conjecture about convergence in the infinite case (Section~\ref{sec:conjecture_main_subsection}). The convergence in these simulations is striking but not perfect, as can be seen in the last column, where the fourth order probabilities (purple) do not perfectly match the supposed point of convergence (black); as far as we are aware, this imperfection is likely due to limitations in sample size and the computer's random generation calculations.

\section{Applications to Single Hypothesis Testing}\label{sec:hypothesis_testing}

We now turn toward the application of our results to single (\ie, null) hypothesis testing. Often, a researcher wishes to test a \emph{single} null hypothesis, such as ``the coin is fair'' or ``the treatment is more than 70\% effective.'' The researcher runs an experiment and wishes to determine whether the result counts as evidence for or against the null hypothesis and, if so, quantify how strong the evidence is. There are presently two major ways to get numeric results for this sort of hypothesis testing: p-values and Bayes factor analysis. That is, when a researcher wants to quantitatively judge whether an experimental result counts as evidence against a null hypothesis, they will employ either some kind of p-value-granting statistical test or obtain a Bayes factor. However, both p-values and standard Bayes factor analysis techniques have serious conceptual issues.

\subsection{Defining Statistical Evidence}

Commonly, researchers do not intend to rely solely on their one experiment to determine their beliefs on a subject. Rather, they intend to use an experimental result as one piece of \emph{evidence} amidst a wider collection of other pieces of evidence -- perhaps some of which will only be acquired in the future.

To our knowledge, the best mathematical, probabilistic interpretation of what it means for a result $R$ to count as \emph{evidence} for or against a null hypothesis (and how strong that evidence is), is how much the result $R$ increases or decreases one's confidence that a hypothesis is true. Expressed mathematically in terms of \emph{precise} probabilities, this means the following: Let $N$ denote the event that the null hypothesis is true. Then an experimental result $R$ should be understood as evidence in the following sense: Let $P(N)$ and $P(\bar{N}) = 1 - P(N)$ be a person's confidences that the null hypothesis is true/false respectively \emph{before} they see the experimental result. Let $P(N | R)$ and $P(\bar{N} | R) = 1 - P(N | R)$ be that person's confidences that the null hypothesis is true/false \emph{after} they see the experimental result. Then $R$ is considered to have been evidence according to the factor by which it changed the person's confidence:

$$\text{evidence for/against }N = \left ( \frac{P(N | R)}{P\left(\bar{N} | R\right)} \middle /  \frac{P(N)}{P\left(\bar{N}\right)} \right ) = \frac{P(R | N)}{P\left(R | \bar{N}\right)}$$

This is to say, if we are going to try to treat evidence as a matter of mathematics and probabilities, then the evidence for/against $N$ is best understood as something equivalent to the Likelihood Ratio, $\frac{P(R | N)}{P\left(R | \bar{N}\right)}$, because the likelihood ratio quantifies how much the result $R$ would change a person's numeric confidence in $N$. This viewpoint is cogently and thoroughly argued by Royal~\citeyearpar{royall1997statistical}. For the purposes of the present work it should be sufficient to note that our stated reason for using likelihood ratios as a valid numeric quantification of evidence: the likelihood ratio quantifies how much a person's numeric confidence would change given the experimental result.

\subsection{If Likelihood Ratios are the Answer, What is the Rest of this Paper For?}

We, along with others, believe that likelihood ratios are the best mathematical formulation of evidence in a probabilistic context \emph{when working with precise probabilities}. \emph{\textbf{However}}, as far as we are aware, there is no good way to acquire a likelihood ratio concerning only a single (null) hypothesis. Put another way, there is no good way to measure evidence about a single hypothesis when using precise probabilities. This is because a likelihood ratio in single-hypothesis testing requires both a probability of the result $R$ given the null hypothesis ($P(R | N)$) \emph{and} a probability of the result $R$ given that the null hypothesis is false ($P\left(R | \bar{N}\right)$). Usually, however, the alternative hypothesis $\bar{N}$ is too vague (\ie, imprecise) to provide this latter probability. For example, consider the question, ``What is the probability of flipping a coin ten times and getting two heads given that the coin is \emph{not} a fair 50-50 coin?'' This question has no answer; $P\left(\cdot | \bar{N}\right)$ is undefined.

Typical Bayes factor analysis purports to provide an answer via its use of prior probabilities over the parameters of the coin in the event that the coin is not 50-50. This technique amounts to a providing a likelihood ratio, but it obtains the likelihood ratio by using a prior over the space of alternative hypotheses to the null hypothesis. We and many others find this technique suspect; in particular, we have found previous justifications for the \emph{choice} of prior to be undermotivated. However, \emph{if} the choice of prior can be sufficiently justified, then Bayes factors seem like a valid measure of evidence, and the prior is effectively a privileged second order confidence function.

\subsection{Using Imprecise Probabilities}

One might think that with imprecise probabilities, we can simply look at the collection of likelihood ratios obtained from each individual probability function in our first order credal set. However, this also is typically very inconclusive, as some probability functions in our first order credal set will consider a test result evidence in favor of the null hypothesis, while others will consider it evidence against the null hypothesis.

Additionally, perhaps our first order event space might not even include the null hypothesis as an \emph{event} to have evidence for or against. Rather, the null hypothesis would simply be a single probability function within our first order credal set, and it would be our second order credal sets that would have our null hypothesis in their event space. Regardless, in this case we still would get that our second order credal set produces likelihood ratios both for and against the null hypothesis.

Rather than looking at an agent's first or second order confidence updates to obtain a measure of evidence, we propose looking at an agent's \textit{highest} order confidence updates. After all, mathematical models of evidence are supposed to model the way an \emph{idealized} agent updates their confidences, and an idealized agent's highest order confidences represent their all-things-considered confidences. Thus, we suggest it makes the most sense to consider an idealized agent's $i^\textrm{th}$ order confidences in the limit as $i \rightarrow \infty$. If those idealized confidences show a particular shift after the agent sees an experimental result, then a quantification of that shift represents the direction and strength of the evidence for that idealized agent.

\subsection{The Higher Order Alternative}\label{sec:our_testing_method}

As just discussed above, we choose to model the strength of statistical evidence as the shift that occurs in an idealized agent's highest (\ie, limiting) order of confidence.

In particular, we specify that the agent sets their first order credal set to contain \emph{all} the relevant hypotheses the experimenter could consider, and then the agent sets all their higher order confidences to be completely agnostic about their first order credal set. This scenario means the agent is considering the experiment, the hypotheses, and nothing more. Importantly, this scenario coincides exactly with our conjecture from Section~\ref{sec:conjecture_main_subsection}, which states that an agent's highest order confidences converge to be identical to a uniform distribution over the agent's first order credal set. In particular, this uniform distribution is defined with respect to a specific distance metric over the space of probability functions rather than a specific numeric parametrization.

Thus, in the end, our method comes out to be functionally equivalent to using Bayes factors with $U_\text{TV}^2$ as the prior. However, \emph{our choice of prior} $U_\text{TV}^2$ is obtained via higher order credal convergence, and is defined implicitly via the statistical distance metric, total variation. Furthermore, in our setup, $U_\text{TV}^2$ does not represent the agent's only prior per-say; even if our conjecture is true, some of the agent's highest order confidence functions will disagree with $U_\text{TV}^2$. Rather than representing an idealized agent's single, prior, $U_\text{TV}^2$ is the prior which represents the consensus or bulk of the idealized agent's confidences.

Importantly, we do not think that our setup for measuring evidence is \emph{perfect} (see the objection in Section~\ref{sec:ideal_vs_us}); rather, we think that our methodology is the best available -- the \textit{least bad}.

We call our measure of evidence the Highest Order Credal Shift Ratio (\ie, the HOCS ratio). For a first order probability function $P^1$, the highest order credal shift from observing an event $E \in \mathcal{F}^1$ can be quantified via:

$$\text{HOCS ratio} = \frac{U_\text{TV}^2(E\ |\ \{P^1\})}{U_\text{TV}^2(E\ |\ \mathcal{C}^1 \setminus \{P^1\})} = \frac{P^1(E)}{U_\text{TV}^2(E\ |\ \mathcal{C}^1 \setminus \{P^1\})}$$

We now proceed to provide an example of our model in action, assuming the truth of our conjecture.

\subsection{Example}






\subsubsection{Binomial Test}

We shall flip a coin $10$ times and assume that the coin tosses are independent. We can now consider the space of probability functions over the result, all the binomial distributions for $n = 10$. We set all higher order credal sets to be completely agnostic.

Given our conjecture, we can express almost all the highest order confidences via the second order uniform distribution $U_\text{TV}^2$. We already depicted the distribution that $U_\text{TV}^2$ expresses over the binomial parameter $p$ and the distribution this implies over heads in Figure~\ref{fig:U_TV_2_binomial}. However, how would one go about calculating it and obtaining that figure?

The full details on deriving $U_\text{TV}^2$ are in Section~\ref{sec:uniformity}. However, we outline the process as follows:

\begin{enumerate}
    \item Obtain a simple parametrization of our hypothesis space; see Section~\ref{sec:simple_parametrizations} for the definition, as not all parametrizations qualify as simple. The definition of a simple paramatrization mostly boils down to the following principle: ``If $x$ is a valid value for parameter $p$, then $x$ must always be a valid value for parameter $p$ regardless of what the other parameters are set to.'' In other words, the parameter space must be an orthotope (a $k$-dimensional rectangle).
    
    In this example, the hypothesis space is the set of binomials, and we can obtain a simple parametrization via the proportion parameter $p$.
    \item Define $U_\text{TV}^2$'s probability density function (PDF) over our hypothesis space. In this case, we have only one parameter, $p$, so the PDF is proportional to:
    $$PDF(p) \propto \lim_{q \rightarrow p} \frac{\text{TV}\left(\text{Binom}(10, p), \text{Binom}(10, q)\right)}{|q - p|}$$
    If the set of hypotheses had a simple parametrization with more than one parameter, then we would need to make the PDF proportional to the product of the above kind of limit for each different parameter.
\end{enumerate}

We can now obtain HOCS ratios by using this PDF as if it were a prior for Bayes factors. What sort of evidence do we get in this scenario?

For any $P^1 \in \mathcal{C}^1$, we have that $U_\text{TV}^2(\{P^1\}) = 0$, so we get that $U_\text{TV}^2(\text{Result}\ |\ \mathcal{C}^1 \setminus \{P^1\}) = U_\text{TV}^2(\text{Result})$. This means our quantification of evidence will be:

$$\text{Evidence about }P^1 = \frac{P^1(\text{Result})}{U_\text{TV}^2(\text{Result})}$$

When the HOCS ratio is greater than 1, the prior says that the result is evidence in favor of $P^1$; a HOCS ratio of less than 1 is evidence against $P^1$.

\begin{figure}[t]
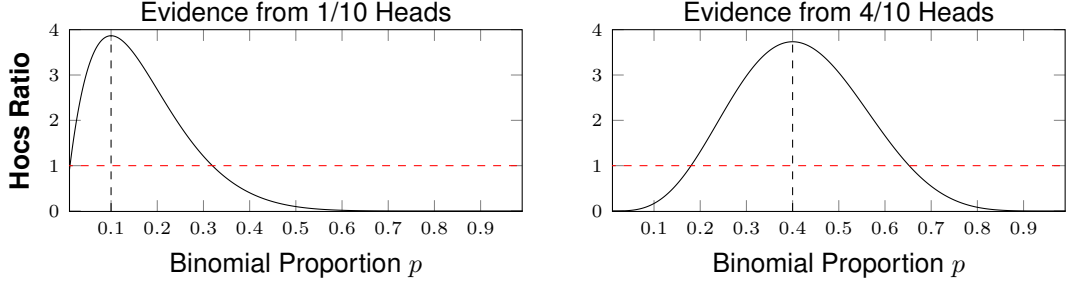

\centering
\input{data/binom_hocs_data_1_10}
\input{data/binom_hocs_data_4_10}

\begin{tikzpicture}
\begin{groupplot}[
        group style={
            group size=2 by 1, 
            group name=hocsplot,
            horizontal sep=1.2cm,
            vertical sep=.5cm,
            xticklabels at=edge bottom,
            xlabels at=edge bottom,
        },
    ylabel style={align=center},
    height=4.0cm,
    width=7.6cm,
    xmin=0.01,
    xmax=0.99,
    ymin=-0.01,
    ymax=4.0,
    tick label style={font=\scriptsize}, 
    every axis title/.style={at={(.5,1.1)},font=\footnotesize},
    xlabel={\footnotesize{\textsf{Binomial Proportion $p$}}},
    xtick={0, .1, .2, .3, .4, .5, .6, .7, .8, .9, 1},
    xmajorticks=true,
]

   \nextgroupplot[
       title={\textsf{Evidence from 1/10 Heads}},
       legend to name={BinomHocsLegend},
       ylabel={\footnotesize{\textsf{~}} \\ \footnotesize{\textbf{\textsf{Hocs Ratio}}}},
   ]
   \addplot[color=black] table [
   x=proportion, y=hocsratio
   ] {\hocsdataoneoften};\label{plots:hocs_1_10}
   
\addplot[dashed,color=red] coordinates {
 (0, 1)
 (1, 1)
};\label{plots:hocs_1_10_center}

\addplot[dashed,color=black] coordinates {
 (0.1, 0)
 (0.1, 3.866381637)
};\label{plots:hocs_1_10_peak}

\nextgroupplot[
       title={\textsf{Evidence from 4/10 Heads}},
       ylabel=\empty, 
   ]
   \addplot[color=black] table [
   x=proportion, y=hocsratio
   ] {\hocsdatafouroften};\label{plots:hocs_4_10}
   
\addplot[dashed,color=red] coordinates {
 (0, 1)
 (1, 1)
};\label{plots:hocs_4_10_center}

\addplot[dashed,color=black] coordinates {
 (0.4, 0)
 (0.4, 3.733746474)
};\label{plots:hocs_4_10_peak}

\end{groupplot}

\end{tikzpicture}
\caption{Highest Order Credal Shift Ratios (HOCS ratios) for two different coin tossing tests results. We have the result for 1 heads out of 10 tosses on the left and 4 of 10 on the right. The HOCS ratios are likelihood functions normalized by the second order uniform distribution, $U_\text{TV}^2$. If our highest-order confidence convergence conjecture is true (Section~\ref{sec:conjecture}), the HOCS ratios quantify the amount by which an idealized agent would update their highest order confidences toward/away-from the possible heads/tails proportions.}\label{fig:binomial_hocs}
\end{figure}

We consider two different results: getting 1 heads and getting 4 heads. We depict the way this provides evidence in Fig.~\ref{fig:binomial_hocs}. These plots are almost identical to plots of likelihood functions; in fact, Fig.~\ref{fig:binomial_hocs} is just the likelihood function scaled via dividing by $U_\text{TV}^2(\text{Result})$. However, while a plot of a likelihood function shows how an experimental result favors one specific hypothesis vis-a-vis another specific hypothesis, it does not provide a determination of evidence for or against a hypothesis simplicater. In Fig.~\ref{fig:binomial_hocs}, the threshold of 1 indicates which hypotheses (\ie, which probability functions) have been supported in the eyes of almost the entirety of our idealized agent's highest order confidence functions.

\section{Objections}\label{sec:objections}

In this section we outline various objections to our ideas and offer some initial responses.

\subsection{Is TV-Uniform Truly Privileged? What About Other Distance Metrics?}

We have asserted that the (conjectured) higher order convergence to $U_\text{TV}^2$ is a natural consequence of our model, and thus $U_\text{TV}^2$ has a privileged, \ie, natural, status. However, we \emph{chose} to use the total variation (TV) distance metric. What if we chose another metric? Would we get a different uniform? For example, if we chose the Hellinger distance metric $H$, would we get $U_H^2$ as the natural uniform? (Note: See Shemyakin for a discussion of such Hellinger priors~\citeyearpar{shemyakin2014hellinger}.). Also, the Hellinger and Total Variation metrics form the same topology on the space of probability functions, so how do we know that our model converges to $U_\text{TV}^2$ and not $U_H^2$? If we chose Fisher Information, perhaps we would have converged on the Jeffereys prior. Not to mention, what about Jensen and Shannon's information theoretic distance? Etc.

Our reply to these concerns has one main component: Trying to define a higher order uniform like we do in Section~\ref{sec:uniformity} is often \emph{impossible} with metrics like Hellinger and Jensen Shannon. This is because the change in metric relative to the change in parameter is \emph{infinite} when moving away from a distribution which places all probability on a singleton event. For example, if we considered a yes/no event and all probability functions over those two possibilities, then the first order probability function which assigns probability 1 to yes would get unbounded probability when trying to define a second order uniform with the Hellinger metric or Jensen Shannon distance. As far as we understand, the Jeffereys prior has the same issue; when people wish to use a Jeffereys prior they must exclude the extreme hypotheses which assign all probability to a singleton event. Even then, Jeffereys and these other uniforms are improper priors, meaning they are not technically probability distributions.

However, we note that one valid alternative to our choice of total variation \emph{might} be a generalization of Euclidean distance via a square root of an $f$-divergence. The Euclidean metric, rather than Manhattan distance, is crucial for our notions of uniformity in 3D space, and it can be generalized via an $f$-divergence to accommodate probability functions over an infinite space. We are unsure as of now whether such a distance measure would be a distance \emph{metric}, but if it is, it may be a reasonable alternative to TV. We find that in finite cases Euclidean distance produces very similar results to TV.

\subsection{Is Evidence for a Theoretical, Completely Agnostic Agent also Evidence for Real Humans?}\label{sec:ideal_vs_us}

Part of the goal of this paper was to say whether an experimental result was evidence for or against a null hypothesis simplicater, which is equivalent to saying whether an experimental result is evidence for or against a null hypothesis ($P_N$) vis-a-vis the alternative hypothesis, where the alternative hypothesis simply states that the null hypothesis is false ($\mathcal{C}^1 \setminus \{P_N\}$). As has been believed for a very long time, such a question cannot be answered validly without reference to some kind of prior over the first order credal set. In this paper, we (via a conjecture) provided a distribution which can act as a privileged prior because it represents almost all of an idealized agent's highest order confidences. However, while we argue that this prior is privileged, we do not argue that it is uninformative.

Thus, while perhaps a truly idealized agent with \emph{no other} information about the experiment would have confidence convergence to second order total variation uniform ($U_\text{TV}^2$), we humans \emph{do} have other information from our life experience, and thus perhaps we should not express evidence in terms of $U_\text{TV}^2$. Perhaps we should not believe using $U_\text{TV}^2$ as a prior describes what evidence an experiment gives \emph{to us}.

As of now, we simply accept this critique. However, we point out that a uniquely privileged prior still offers a uniquely privileged notion of what evidence an experiment provides in and of itself. Furthermore, if $U_\text{TV}^2$ is truly a uniquely privileged prior for an agent with no previous experiences, then $U_\text{TV}^2$, coupled with an agent's experiences expressed as an event $E$, provides a uniquely privileged prior for that agent: $m_E^2(U_\text{TV}^2)$.

\section{Conclusion}

We provided a formal model to generalize confidence updates with credal sets by allowing for higher order confidences via higher order credal sets. This model shows promise both toward addressing standard issues with credal set models and toward providing an acceptable version of single hypothesis testing (if such a thing is even possible). If our conjecture is true, it provides us with a prioritized, natural prior function over a space of hypotheses -- namely, uniformity with respect to the total variation metric.





\end{document}